\input pictex.tex
\input amstex.tex
\documentstyle{amsppt}
\magnification=1200
 \pagewidth{16.3truecm}
 \pageheight{24truecm}
 \nologo

\NoBlackBoxes
\def\div{\operatorname{div}}
\def\mod{\operatorname{mod}}

\def\supp{\operatorname{supp}}

\def\gcd{\operatorname{gcd}}
\def\Pic{\operatorname{Pic}}

\refstyle{A} \widestnumber\key{ACLM2}
 \topmatter
\title
Monodromy eigenvalues and zeta functions with differential forms
\endtitle
\author
Willem Veys  
\endauthor
\address K.U.Leuven, Departement Wiskunde, Celestijnenlaan 200B,
         B--3001 Leuven, Belgium  \endaddress
\email wim.veys\@wis.kuleuven.be  \newline
 http://www.wis.kuleuven.be/algebra/veys.htm
\endemail
 \abstract
For a complex polynomial or analytic function $f$, there is a
strong correspondence between poles of the so-called local zeta
functions or  complex powers $\int |f|^{2s}\omega$, where the
$\omega$ are $C^{\infty}$ differential forms with compact support,
and eigenvalues of the local monodromy of $f$. In particular
Barlet showed that each monodromy eigenvalue of $f$ is of the form
$\exp(2 \pi \sqrt{-1} s_0)$, where $s_0$ is such a pole. We prove
an analogous result for similar $p$-adic complex powers, called
Igusa (local) zeta functions, but mainly for the related
algebro-geometric topological and motivic zeta functions.
\endabstract
 \endtopmatter

\document

\noindent {\bf Introduction}
\bigskip
\noindent {\bf 0.1.} Let $f : X \rightarrow \Bbb C$ be a
non-constant analytic function on an open part $X$ of $\Bbb C^n$.
We consider $C^\infty$ functions $\varphi$ with compact support on
$X$ and the corresponding differential forms $\omega = \varphi dx
\wedge d \bar x$. Here and further $x = (x_1,\cdots,x_n)$ and $dx
= dx_1 \wedge \cdots \wedge dx_n$. For such $\omega$ the integral
$$Z(f,\omega;s) := \int_X | f(x) |^{2s} \omega,$$
where $s \in \Bbb C$ with $\Re (s) > 0$, has been the object of
intensive study. One verifies that $Z(f,\omega;s)$ is holomorphic
in $s$. Either by resolution of singularities [At], [BG], or by
the theory of Bernstein polynomials [Be], one can show that it
admits a meromorphic continuation to $\Bbb C$, and that all its
poles are among the translates by $\Bbb Z_{< 0}$ of a finite
number of rational numbers. Combining results of Barlet [Ba2],
Kashiwara [Ka] and Malgrange [Ma2], the poles of (the extended)
$Z(f,\omega;s)$ are strongly linked to the eigenvalues of (local)
monodromy at points of $\{ f = 0 \}$; see \S 1 for the concept of
monodromy.

\bigskip
\proclaim{Theorem} (1) If $s_0$ is a pole of $Z(f,\omega;s)$ for
some diffential form $\omega$, then $\exp(2 \pi \sqrt{-1} s_0)$ is
a monodromy eigenvalue of $f$ at some point of $\{ f = 0 \}$.

(2) If $\lambda$ is a monodromy eigenvalue of $f$ at a point of
$\{ f = 0 \}$, then there exists a differential form $\omega$ and
a pole $s_0$ of $Z(f,\omega;s)$ such that $\lambda = \exp(2 \pi
\sqrt{-1} s_0)$.
\endproclaim

There are also more precise local versions in a neighbourhood of a
point of $\{f=0\}$.
 Similar results hold for a real
analytic function $f : X(\subset \Bbb R^n) \rightarrow \Bbb R$ and
integrals $\int_{X \cap \{ f > 0 \}} f^s \varphi dx$; we refer to
e.g. [Ba3], [Ba4], [Ba6], [BM], [JM].

\bigskip
\noindent {\bf 0.2.} Let now $f : X \rightarrow \Bbb Q_p$ be a
non-constant ($\Bbb Q_p$-)analytic function on a compact open $X
\subset \Bbb Q^n_p$, where $\Bbb Q_p$ denotes the field of
$p$-adic numbers. Let $| \cdot |_p$ and $| dx |$ denote the
$p$-adic norm and the Haar measure on $\Bbb Q^n_p$, normalized in
the standard way. The $p$-adic integral
$$Z_p(f;s) := \int_X | f(x)|^s_p |dx|,$$
again defined for $s \in \Bbb C$ with $\Re(s) > 0$, is called the
($p$-adic) Igusa zeta function of $f$. Using resolution of
singularities Igusa [Ig1] showed that it is a rational function of
$p^{-s}$; hence it also admits a meromorphic continuation to $\Bbb
C$. In this context there is an intriguing conjecture of Igusa
relating poles of (the extended) $Z_p(f;s)$ to eigenvalues of
monodromy. More precisely, let $f$ be a polynomial in $n$
variables over $\Bbb Q$. Then we can consider $Z_p(f;s)$ for all
prime numbers $p$ (taking $X = \Bbb Z^n_p)$.

\bigskip
\noindent \proclaim{Monodromy conjecture {\rm [De1]}} For all
except a finite number of $p$, we have that, if $s_0$ is a pole of
$Z_p(f;s)$, then $\exp(2 \pi \sqrt{-1} s_0)$ is a monodromy
eigenvalue of $f : \Bbb C^n \rightarrow \Bbb C$ at a point of $\{
f = 0 \}$.
\endproclaim

\bigskip
\noindent This conjecture was proved for $n=2$ by Loeser [Lo2].
There are by now various other partial results [ACLM1], [ACLM2],
[Lo3], [RV], [Ve1], [Ve6]. (We took $\Bbb Q_p$ for simplicity of
notation; everything can be done over finite extensions of $\Bbb
Q_p$.)

\bigskip
\noindent{\bf 0.3.} There are various `algebro-geometric' zeta
functions, related to the $p$-adic Igusa zeta functions: the
motivic, Hodge and topological zeta functions, for which we refer
to \S 1. Here we just mention that the motivic zeta function
specializes to the various $p$-adic Igusa zeta functions (for
almost all $p$). For those zeta functions a similar monodromy
conjecture can be stated; and analogous partial results are valid.

\bigskip
\noindent{\bf 0.4.} We should note that for the complex (and real)
integrals in (0.1) there are more precise results of Bernstein and
Barlet, involving roots of the Bernstein polynomial of $f$
(instead of monodromy eigenvalues). Similarly there is a finer
conjecture for the poles of Igusa and related zeta functions,
relating them to roots of the Bernstein polynomial [De1], [DL1].
However, the results of this paper do not involve Bernstein
polynomials, so we just refer the interested reader to [Be],
[Ba1], [Ba2], [Ba3], [Ig2], [Lo1], [Lo2], [Lo3].

\bigskip
\noindent{\bf 0.5.} As in the complex (or real) case, one
associates $p$-adic Igusa zeta functions, and also motivic, Hodge
and topological zeta functions, to a function $f$ {\it and} a
differential form $\omega$. In this `algebro-geometric' context
one considers algebraic differential forms $\omega$; see \S 1.

To our knowledge a possible analogue of Theorem 0.1(2) in the
context of $p$-adic and the related `algebro-geometric' zeta
functions was not studied before in the literature. For instance
let $f$ be a polynomial over $\Bbb Q$ satisfying $f(0) = 0$, and
let $\lambda$ be a monodromy eigenvalue of $f$ at $0$. Does there
exist a compact open neighbourhood $X$ of $0$ and an algebraic
differential form $\omega$ such that (the meromorphic continuation
of) $\int_X |f(x)|^s_p |\omega|_p$ has a pole $s_0$ satisfying
$\lambda = \exp (2 \pi \sqrt{-1} s_0)$? (If $\omega = g(x)dx$ for
some polynomial $g$ over $\Bbb Q$, the integral above is just
$\int_X |f(x)|^s_p |g(x)|_p |dx|$.)

\bigskip
\noindent {\bf 0.6.} We will concentrate in this paper on the
analogous question for the topological zeta function, since a
positive answer in this context automatically yields a positive
answer in the context of Hodge and motivic zeta functions
 (see \S1), and also for Igusa zeta functions (by [DL1, Th\'eor\`eme
2.2]).
 We show for instance (Theorem 3.6):

\bigskip
\noindent \proclaim{Theorem} Let $f : (\Bbb C^n,0) \rightarrow
(\Bbb C,0)$ be a nonzero polynomial function (germ). Let $\lambda$
be a monodromy eigenvalue of $f$ at $0$. Then there exists a
differential $n$-form $\omega$ and a point $P \in \{ f = 0 \}$,
close to $0$, such that the (local) topological zeta function at
$P$, associated to $f$ and $\omega$, has a pole $s_0$ satisfying
$\exp(2 \pi \sqrt{-1} s_0) = \lambda$.

If $f^{-1}\{0\}$ has an isolated singularity at $0$, then we can
take $0$ itself as point $P$.
\endproclaim

\bigskip
For $n=2$, we construct such $\omega$ in $\S$2 using so-called
curvettes. In arbitrary dimension we follow a similar approach,
for which we first introduce a higher dimensional version of this
notion (Proposition 3.2).

\bigskip
\noindent {\bf 0.7.} The zeta functions associated to $f$ and the
constructed $\omega$ in the theorem above can have {\sl other}
poles that don't induce monodromy eigenvalues of $f$. So for those
zeta functions the analogue of Theorem 0.1(1) is (unfortunately)
not true. It would be really interesting to have a complete
analogue of Theorem 0.1, roughly saying that the monodromy
eigenvalues of $f$ correspond precisely to the poles of the zeta
functions associated to $f$ and some finite list of differential
forms $\omega$ (including $dx$). Of course this would be a lot
stronger than the (in arbitrary dimension) still wide open
monodromy conjecture.

However, we indicate some examples where such a correspondence
holds, for instance $f = y^a - x^b$ with $\gcd(a,b) = 1$.

\bigskip
\noindent {\sl Acknowledgement.} We would like to thank the
referee for an interesting suggestion.

\bigskip
\bigskip
\noindent {\bf 1. Monodromy and zeta functions}
 \bigskip
  \noindent
{\bf 1.1.} Let $f : \Bbb C^n \rightarrow \Bbb C$ be a non-constant
polynomial function satisfying $f(b) = 0$. Let $B \subset \Bbb
C^n$ be a small enough ball with centre $b$; the restriction
$f|_B$ is a topological fibration over a small enough pointed disc
$D \subset \Bbb C \backslash \{ 0 \}$ with centre $0$. The fibre
$F_b$ of this fibration is called the (local) Milnor fibre of $f$
at $b$; see e.g. [Mi]. The counterclockwise generator of the
fundamental group of $D$ induces an automorphism of the
cohomologies $H^q(F_b, \Bbb C)$, which is called the (local)
monodromy of $f$ at $b$. By a monodromy eigenvalue of $f$ at $b$
we mean an eigenvalue of the monodromy action on a least one of
the $H^q(F_b, \Bbb C)$. It is well known that $H^q(F_b,\Bbb C) =
0$ for $q \geqslant n$, and that all monodromy eigenvalues are
roots of unity.

Let $P_q(t)$ denote the characteristic polynomial of the monodromy
action on $H^q(F_b,\Bbb C)$. If $f = \prod_j f^{N_j}_j$ is the
decomposition of $f$ in irreducible components and $d := \gcd_j
N_j$, then $P_0(t) = t^d - 1$.

When $b$ is an isolated singularity of $f^{-1} \{ 0 \}$, then
$H^q(F_b,\Bbb C) = 0$ for $q \ne 0, n-1$; and $P_0(t) = t-1$.

\bigskip
\noindent {\bf 1.2. Definition.} The {\it monodromy zeta function}
$\zeta_{f,b}(t)$ of $f$ at $b$ is the alternating product of all
characteristic polynomials $P_q(t)$ :
$$\zeta_{f,b}(t) :=  \prod^{n-1}_{q=0} P_q (t)^{{(-1)}^q}.$$
Note that there are also other conventions, see for example [A'C],
[AVG].

\bigskip
\noindent In particular for an isolated singularity the knowledge
of $\zeta_{f,b}(t)$ and of $P_{n-1}(t)$ are equivalent.

\bigskip
\noindent {\bf 1.3.} We recall the following interesting and
useful result, which is maybe not generally known.

\bigskip
\proclaim{Lemma {\rm [De2, Lemma 4.6]}} Let $f : \Bbb C^n
\rightarrow \Bbb C$ be a non-constant polynomial function. If
$\lambda$ is a monodromy eigenvalue of $f$ at $b \in f^{-1} \{ 0
\}$, then there exists $P \in f^{-1} \{ 0 \}$ (arbitrarily close
to $b$) such that $\lambda$ is a zero or pole of the monodromy
zeta function of $f$ at $P$.
\endproclaim

\bigskip
\noindent

It is convenient to recall also the proof in order to see how the
point $P$ is obtained. Let $\Psi_{f,\lambda}$ be the sub-complex
of the complex of nearby cycles of $f$ corresponding to the
eigenvalue $\lambda$, where both are viewed as (shifted) perverse
sheaves. Let $\Sigma$ be the largest analytic set given by the
supports of the cohomology sheaves of $\Psi_{f,\lambda}$. Then, by
perversity of $\Psi_{f,\lambda}$, at a generic point $P$ of
$\Sigma$ the eigenvalue $\lambda$ appears on exactly one
cohomology group of the Milnor fibre of $f$ at $P$.

\bigskip
\noindent {\bf 1.4.} {\sl A'Campo's formula.} Let $f : \Bbb C^n
\rightarrow \Bbb C$ be a non-constant polynomial function
satisfying $f(b) = 0$. Take an embedded resolution $\pi : X
\rightarrow \Bbb C^n$ of $f^{-1} \{ 0 \}$ (that is an isomorphism
outside the inverse image of $f^{-1}\{0\}$). Denote by $E_i, i \in
S$, the irreducible components of the inverse image
$\pi^{-1}(f^{-1} \{ 0 \})$, and by $N_i$ the multiplicity of $E_i$
in the divisor of $\pi^\ast f$. We put $E^\circ_i := E_i \setminus
\cup_{j \ne i} E_j$ for $i \in S$.

\bigskip
\proclaim{Theorem {\rm [A'C]}} Denoting by $\chi(\cdot)$ the
topological Euler characteristic we have
$$\zeta_{f,b}(t) = \prod_i (t^{N_i}-1)^{\chi(E^\circ_i \cap
\pi^{-1} \{b \})}.$$
\endproclaim

\bigskip
\noindent {\bf 1.5.} Another kind of zeta functions are the
topological, Hodge and motivic zeta functions, associated to a
non-constant polynomial function $f : \Bbb C^n \rightarrow \Bbb C$
and a regular differential $n$-form $\omega$ on $\Bbb C^n$. (More
generally one can consider an arbitrary smooth quasi-projective
variety $X_0$ instead of $\Bbb C^n$ and a regular function $f$.)
We will describe these zeta functions in terms of an embedded
resolution of $f^{-1} \{ 0 \} \cup  \div  \omega$. Now we denote
by $E_i, i \in S$, the irreducible components of the inverse image
$\pi^{-1} (f^{-1} \{ 0 \} \cup  \div \omega)$ and by $N_i$ and
$\nu_i - 1$ the multiplicities of $E_i$ in the divisor of
$\pi^\ast f$ and $\pi^\ast \omega$, respectively. We put
$E^\circ_I := (\cap_{i \in I} E_i) \setminus (\cup_{j \notin I}
E_j)$ for $I \subset S$, in particular $E^\circ_\emptyset = X
\setminus (\cup_{j \in S} E_j)$. So the $E^\circ_I$ form a
stratification of $X$ in locally closed subsets.

\bigskip
\noindent {\bf Definition.} The (global) {\sl topological zeta
function} of $(f,\omega)$ and its local version at $b \in \Bbb
C^n$ are
$$Z_{top} (f, \omega; s) := \sum_{I \subset S} \chi (E^\circ_I)
\prod_{i \in I} \frac{1}{\nu_i + sN_i},$$ and
$$Z_{top,b} (f, \omega; s) := \sum_{I \subset S} \chi (E^\circ_I \cap \pi^{-1} \{ b \})
\prod_{i \in I} \frac{1}{\nu_i + sN_i},$$ respectively, where $s$
is a variable.

\bigskip
These invariants were introduced by Denef and Loeser in [DL1] for
`trivial $\omega$', i.e. for $\omega = dx_1 \wedge \cdots \wedge
dx_n$. Their original proof that these expressions do not depend
on the chosen resolution is by describing them as a kind of limit
of $p$-adic Igusa zeta functions. Later they obtained them as a
specialization of the intrinsically defined motivic zeta functions
[DL2]. Another technique is applying the Weak Factorization
Theorem [AKMW], [W\l] to compare two different resolutions. For
arbitrary $\omega$ one can proceed analogously.

It is natural and useful to study these invariants incorporating
such a more general $\omega$, see for example [ACLM1], [ACLM2],
[Ve4]. Note however that there one restricts to the situation
where $\supp(\div \omega) \subset f^{-1} \{ 0 \}$. In the original
context of $p$-adic Igusa zeta functions, see e.g. [Lo2, III 3.5].

\bigskip
\noindent {\bf 1.6.} There are finer variants of these zeta
functions using, instead of Euler characteristics, Hodge
polynomials or classes in the Grothendieck ring of varieties. We
mention for instance the {\it Hodge zeta function}
$$Z_{Hod}(f,\omega;T) = \sum_{I \subset S} H(E^\circ_I; u,v)
\prod_{i \in I} \frac{(uv-1)T^{N_i}}{(uv)^{\nu_i} - T^{N_i}} \in
\Bbb Q (u,v)(T),$$ where $H(\cdot;u,v) \in \Bbb Z[u,v]$ denotes
the Hodge polynomial. Concerning Hodge and motivic zeta functions
we refer to e.g. [DL2], [Ro2], [Ve5] for versions with `trivial
$\omega$'; and to [ACLM1], [ACLM2], [Ve4] involving more general
$\omega$. We just mention that, in contrast with topological zeta
functions, Hodge and motivic zeta functions can be defined
intrinsically as formal power series (in $T$) with coefficients
determined by the behaviour of the arcs on $\Bbb C^n$ with respect
to their intersection with $f^{-1} \{ 0 \}$ and with
$\div(\omega)$. Then one shows that they are rational functions
(in $T$) by proving explicit formulae as above in terms of an
embedded resolution. However, the fact that we allow (and need)
differential forms $\omega$ with $\supp (\div \omega) \not\subset
f^{-1} \{ 0 \}$, causes some technical complications. In order to
avoid these, one can {\it define} Hodge and motivic zeta functions
by a formula as above in terms of an embedded resolution, and use
the Weak Factorization Theorem to show independency of the choice
of resolution. (In fact in this paper we will only use
differential forms $\omega$ for which a given embedded resolution
of $f^{-1} \{ 0 \}$ is also an embedded resolution of $f^{-1} \{ 0
\} \cup \div \omega$.)

\bigskip
\noindent {\bf 1.7.} We now explain why we choose not to give
details here about Hodge and motivic zeta functions. The point is
that, for a given $f$ and $\omega$, the motivic zeta function
specializes to the Hodge zeta function, which in turn specializes
to the topological zeta function. (Note for instance that
$H(\cdot; 1,1) = \chi(\cdot)$.) In particular, a pole of the
topological zeta function will induce a pole of the other two.
(The converse is not clear.) The problem that we want to treat
here is, given a monodromy eigenvalue $\lambda$ of $f$, find a
form $\omega$ such that the zeta function associated to $f$ and
$\omega$ has a pole `inducing $\lambda$'. Therefore in this paper
we focus on the topological zeta function. We will succeed in
proving the desired result for it, implying the analogous result
for the `finer' zeta functions.

\bigskip
\bigskip
\noindent {\bf 2. Curves}

\bigskip
\noindent {\bf 2.1.} We first prove our main result for curves.
Ultimately Theorem 2.4 below will be essentially a special  case
of Theorems 3.5 and 3.6. It is however more precise in the case of
non-isolated singularities. We also believe that it is useful to
treat the case of curves first. The proof is easier and shorter,
and we indicate a fact that is typical for curves.

\bigskip
\noindent {\bf 2.2.} In order to construct appropriate
differential forms $\omega$ we use curve germs (sometimes called
{\it curvettes}) that intersect the exceptional components of an
embedded resolution transversely; we quickly recall this notion.
Let
$$\Bbb A^2 {\overset \pi_1 \to \longleftarrow} X_1 {\overset \pi_2 \to \longleftarrow} X_2
{\overset \pi_3 \to \longleftarrow} ... {\overset \pi_i \to
\longleftarrow} X_i {\overset \pi_{i + 1} \to \longleftarrow}
\cdots {\overset \pi_m \to \longleftarrow} X_m
$$
be a composition $\pi$ of $m$ blowing-ups with $0$ the centre of
$\pi_1$, and the centre of all other $\pi_i$ belonging to the
exceptional locus of $\pi_1 \circ \cdots \circ \pi_{i-1}$. In
other words, all centres are points infinitely near to $0$. Denote
the exceptional curve of $\pi_i$, as well as its strict transform
in $X_m$, by $E_i$. A {\it curvette} $C_i$ of $E_i$ is a smooth
curve (germ) on $X_m$ satisfying $C_i \cdot E_j = \delta_{ij}$ for
all $j = 1, \cdots, m$. So $C_i$ intersects $E_i$ transversely in
a point not belonging to other $E_j$. We denote $\bar C_i :=
\pi(C_i)$ the image curve (germ) of $C_i$ in $(\Bbb A^2,0)$.

We guess that the following should be known.

\bigskip
\proclaim{2.3. Proposition} Let $\pi^\ast \bar C_i = \sum^m_{j=1}
a_{ij} E_j + C_i$ for $i = 1,\cdots,m$. Then the determinant of
the $(m \times m)$-matrix $(a_{ij})$ is equal to 1. In particular
$\gcd_{1 \leq i \leq m} \{ a_{ij} \} = 1$ for all $j$.
\endproclaim

\noindent \demo{Proof} For all $i, \ell \in \{1,\cdots,m\}$ we
have
$$0 = (\pi^\ast \bar C_i) \cdot E_\ell = \sum^m_{j=1} a_{ij} E_j
\cdot E_\ell + \delta_{i \ell}.$$ In other words, the matrix
product $(a_{ij}) \cdot (-E_j \cdot E_\ell)$ is the identity
matrix. Since minus the intersection matrix of the $E_i$ has
determinant $1$, the same is true for $(a_{ij})$. $\square$
\enddemo

\bigskip
\noindent {\sl Note.} It also follows that $(a_{ij})$ is
symmetric, being the inverse of a symmetric matrix.

\bigskip
\proclaim{2.4. Theorem} Let $f : (\Bbb C^2,0) \rightarrow (\Bbb
C,0)$ be a non-zero polynomial function (germ). Let $\lambda$ be a
monodromy eigenvalue of $f$ at $0$, i.e. $\lambda$ is an
eigenvalue for the action of the local monodromy on $H^0(F_0, \Bbb
C)$ or $H^1(F_0,\Bbb C)$. Then there exists a differential 2-form
$\omega$ on $(\Bbb C^2,0)$ such that $Z_{top,0}(f,\omega;s)$ has a
pole $s_0$ satisfying $\exp(2 \pi \sqrt{-1} s_0) = \lambda$.
\endproclaim

\demo{Proof} Take the minimal embedded resolution $\pi : X
\rightarrow (\Bbb C^2,0)$ of the curve (germ) $f^{-1} \{ 0 \}$.
Denote the irreducible components of $\pi^{-1} (f^{-1} \{ 0 \})$,
i.e. the exceptional curves and the components of the strict
transform, by $E_i$ and their multiplicities in the divisors of
$\pi^\ast f$ and $\pi^\ast(dx \wedge dy)$ by $N_i$ and $\nu_i -
1$, respectively. Take for each exceptional curve $E_i$ a
(generic) curvette $C_i \subset X$. Say $\bar C_i := \pi(C_i)$ is
given by the (reduced) equation $g_i = 0$.

\medskip
 We first suppose that $\lambda$ is a pole of the monodromy
zeta function of $f$ at $0$. Say $\lambda$ is a primitive $d$-th
root of unity. By A'Campo's formula there exists an exceptional
curve $E_{j_0}$ with $d|N_{j_0}$ and $\chi(E^\circ_j) < 0$. So
$E_{j_0}$ intersects at least three times other components $E_j$.

We associate to each curvette the following multiplicities:
$$\div(\pi^\ast g_\ell) = \sum_j a_{\ell j} E_j + C_\ell.$$
Note that $a_{\ell j} \ne 0$ only if $E_j$ is exceptional. We take
a differential form $\omega$ of the form $(\prod_\ell
g^{m_\ell}_\ell) dx \wedge dy$. The multiplicity  of $E_i$ in the
divisor of $\pi^\ast \omega$ is $\nu_i - 1 + \sum_\ell a_{\ell i}
m_\ell$. So in particular the candidate pole for
$Z_{top,0}(f,w;s)$ associated to $E_{j_0}$ is
$$s_0 = -\frac{\nu_{j_0} + \sum_\ell a_{\ell j_0}
m_\ell}{N_{j_0}}.$$ We will find suitable $m_\ell$ such that (1)
$\exp(2 \pi \sqrt{-1} s_0) = \lambda$, and (2) $s_0$ is really a
pole.

\bigskip (1) Say $\lambda = \exp (-2 \pi \sqrt{-1} \frac bd)$ with
$1 \leq b \leq d$ and $\gcd(b,d) = 1$. Since $\gcd_\ell \{a_{\ell
j_0} \} = 1$ by Proposition 2.3, there exist $m_\ell \in \Bbb Z$
such that $\nu_{j_0} + \sum_\ell a_{\ell j_0} m_\ell = b
\frac{N_{j_0}}{d}$. Consequently there exist $m_\ell \in \Bbb
Z_{\geq 0}$ satisfying
$$\nu_{j_0} + \sum_\ell a_{\ell j_0}m_\ell = b \frac{N_{j_0}}{d}\, \mod
N_{j_0} ,$$ and we can choose such $m_\ell$ freely in their
congruence class $\mod N_{j_0}$. For all those $m_\ell$ clearly
$\exp(2 \pi \sqrt{-1} s_0) = \exp (-2 \pi \sqrt{-1} \frac bd) =
\lambda$.

\bigskip
(2) The candidate pole for $Z_{top,0}(f,\omega;s)$ associated to a
component $E_i$ of the strict transform is $-\frac1{N_i}$, and is
thus different from $s_0$ for \lq most\rq\ $m_\ell$. The candidate
pole associated to another exceptional $E_i$ is  $- \frac{\nu_i +
\sum_\ell a_{\ell i} m_\ell}{N_i}$. Suppose that it is equal to
$s_0$. Then
$$\frac{\nu_i}{N_i} - \frac{\nu_{j_0}}{N_{j_0}} + \sum_\ell
(\frac{a_{\ell i}}{N_i} - \frac{a_{\ell j_0}}{N_{j_0}})m_\ell = 0.
\tag {$*$}$$ We know that det$(a_{ij}) \ne 0$ by Proposition 2.3;
in particular the vectors $(a_{\ell i})_\ell$ and $(a_{\ell
j_0})_\ell$ cannot be dependent. So at least one of the
coefficients of the $m_\ell$ in ($*$) is non-zero, i.e. ($*$) is
never an empty condition. Consequently `most' sets $(m_\ell)_\ell$
in our allowed lattice satisfy $s_0 \ne - \frac{\nu_i + \sum_\ell
a_{\ell i} m_\ell}{N_i}$ for all $i \ne j_0$.

The residue of $s_0$ (as candidate pole of order 1 for $Z_{top,0}
(f,\omega;s))$ is
$$\frac{1}{N_{j_0}} (\chi (E^\circ_{j_0}) - 1 + \frac{1}{1 +
m_{j_0}} + \sum_i \frac{1}{\alpha_i})$$ where $\alpha_i := \nu_i -
\frac {\nu_{j_0}}{N_{j_0}} N_i + \sum_\ell (a_{\ell i} - a_{\ell
j_0} \frac{N_i}{N_{j_0}}) m_\ell$ for an $E_i$ intersecting
$E_{j_0}$. (See Figure 1.)

Since $\chi(E^\circ_{j_0}) - 1 \ne 0$, this expression is never
identically zero as function in the $m_\ell$, and hence nonzero
for `most' choices of $(m_\ell)_\ell$.

\bigskip
Secondly, if an eigenvalue $\lambda$ of $f$ at $0$ is not a pole
of the monodromy zeta function, it must be an eigenvalue on
$H^0(F_0,\Bbb C)$. By (1.1) the order $d$ of $\lambda$ as root of
unity must divide all $N_j$ associated to components of the strict
transform. But then $d$ divides all $N_i$.

Pick now any exceptional $E_{j_0}$ with $\chi (E^\circ_{j_0}) < 0$
and proceed as in the previous case to construct a suitable
$\omega$ and a pole $s_0$ of $Z_{top,0}(f, \omega; s)$ with
$\exp(2 \pi \sqrt{-1} s_0) = \lambda$. (There is always an
exceptional $E_i$ with $\chi(E^\circ_i) < 0$, except in the
trivial case where $f^{-1} \{ 0 \}$ is smooth or has normal
crossings at $0$. But then also the theorem is quite trivial, see
Example 2.6.) $\square$
\enddemo

\vskip 1truecm
 \centerline{
\beginpicture
\setcoordinatesystem units <.5truecm,.5truecm>
 \putrule from 0 1 to 8.5 1
 \putrule from 1.5 0 to 1.5 4
 \putrule from 4 0 to 4 4
 \putrule from 5 0 to 5 4
 \putrule from 7 0 to 7 4
 \put {$\dots$} at 6 3
  \put {$C_{j_0}$} at .8 4
  \put {$E_{j_0}$} at 9.3 1
  \put {$E_i$} at 6 4
\endpicture}
   \vskip 1truecm
  \centerline{\smc Figure 1}
 \vskip 1truecm

\bigskip
\noindent {\bf 2.5.} {\sl Note.} The zeta functions $Z_{top,0}
(f,\omega;s)$ constructed in the proof above can in general have
other poles that {\sl don't} induce monodromy eigenvalues of $f$.
It would be interesting to investigate the validity of the
following statement, or its analogue for Hodge and motivic zeta
functions.

{\sl Let $f : (\Bbb C^2,0) \rightarrow (\Bbb C,0)$ be a non-zero
polynomial function (germ). There exist regular differential
$2$-forms $\omega_1,\cdots,\omega_r$ on $(\Bbb C^2,0)$ such that

(1) if $s_0$ is a pole of a zeta function
$Z_{top,0}(f,\omega_i;s)$, then $\exp (2 \pi \sqrt{-1} s_0)$ is a
monodromy eigenvalue of $f$ at $0$, and

(2) for each monodromy eigenvalue $\lambda$ of $f$ at $0$, there
is a differential form $\omega_i$ and a pole $s_0$ of
$Z_{top,0}(f,\omega_i;s)$ such that $\exp(2 \pi \sqrt{-1} s_0) =
\lambda$.}

\bigskip
We present some examples of this principle.

\bigskip
\noindent {\bf 2.6.} {\sl Baby example.} (1) Let $f = x^N$ on
$(\Bbb C^2,0)$. The monodromy eigenvalues of $f$ at $0$ are
$\exp(2 \pi \sqrt{-1} \frac bN)$ with $1 \leq b \leq N$. Take
$\omega_b := x^{b-1} dx \wedge dy$ for $b = 1,\cdots,N$. We have
$Z_{top,0}(f,\omega_b;s) = \frac{1}{b + sN}$ with unique pole $s_0
= - \frac bN$.

(2) Let $f = x^{dN} y^{dN^\prime}$ on $(\Bbb C^2,0)$ with
$\gcd(N,N^\prime) = 1$. The monodromy eigenvalues of $f$ at $0$
are $\exp(2 \pi \sqrt{-1} \frac bd)$ with $1 \leq b \leq d$. Take
$\omega_b := x^{bN-1} y^{bN'-1} dx \wedge dy$ for $b =
1,\cdots,d$. We have
$$Z_{top,0}(f,\omega_b;s) = \frac{1}{(bN + sdN)(bN^\prime +
sdN^\prime)} = \frac{1}{NN^\prime (b + sd)^2}$$ with unique pole
(of order 2) $s_0 = - \frac bd$.

\vskip 1truecm
 \centerline{
\beginpicture
\setcoordinatesystem units <.5truecm,.5truecm>
 \setdashes
 \putrule from -2 0 to 2 0
 \putrule from 0 -2 to 0 2
 \setsolid
\ellipticalarc axes ratio 2:3  70 degrees from 0 0 center at 0 3
\ellipticalarc axes ratio 2:3  -70 degrees from 0 0 center at 0 -3
\put {$\bullet$} at 0 0
 \put {$\bar C_1$} at -.6 -2
 \put {$\bar C_2$} at -2 .6
 \put {$f^{-1}\{0\}$} at 2.9 2.9
\put {$\longleftarrow$} at 6 0
 \put {$\pi$} at 7 .5
\setcoordinatesystem units <.5truecm,.5truecm> point at -19 0
 \putrule from 0 -5 to 0 5
 \putrule from -1 4 to 3 4
 \putrule from 2 5 to 2 1
 \putrule from 1 2 to 4 2
 \put{$\ddots$} at 4.5 1.2
  \putrule from 5 0 to 8 0
 \putrule from 7 1 to 7 -3
  \putrule from 1 -4 to -3 -4
 \putrule from -2 -5 to -2 -1
 \putrule from -1 -2 to -4 -2
 \put{$\ddots$} at -4.5 -.8
  \putrule from -5 0 to -8 0
 \putrule from -7 -1 to -7 3
\setdashes
 \putrule from 6 -2 to 9.2 -2
 \putrule from -6 2 to -9.2 2
 \setlinear \plot -1.8 1  1.8 -1 /
 \setsolid
 \put {$E_1$} at -7.4 3.5
  \put {$C_m$} at -2 1.6
  \put {$C_1$} at -9.2 1.4
  \put {$E_2$} at 7.5 -3.5
   \put {$E_m$} at .8 -2.3
  \put {$C_2$} at 9.4 -1.5
\endpicture}
 \vskip 1truecm
  \centerline{\smc Figure 2}
 \vskip 1truecm

\bigskip
\proclaim{2.7. Proposition} Let $f = y^p - x^q$ on $(\Bbb C^2,0)$
with $2 \leq p < q$ and $\gcd(p,q) = 1$. Take $\omega_{ij} :=
x^{i-1}y^{j-1} dx \wedge dy$ for $1 \leq i \leq q-1$ and $1 \leq j
\leq p-1$.

(1) If $s_0$ is a pole of $Z_{top,0}(f,\omega_{ij};s)$ for some
$\omega_{ij}$, then $\exp(2 \pi \sqrt{-1} s_0)$ is a monodromy
eigenvalue of $f$ at $0$.

(2) If $\lambda$ is a monodromy eigenvalue of $f$ at $0$, then
there is a form $\omega_{ij}$ and a pole $s_0$ of
$Z_{top,0}(f,\omega_{ij};s)$ such that $\exp(2 \pi \sqrt{-1} s_0)
= \lambda$.
\endproclaim

\noindent \demo{Proof} Let $\pi : X \rightarrow (\Bbb C^2,0)$ be
the minimal embedded resolution of $f^{-1} \{ 0 \}$, see Figure 2.
We can consider the strict transforms $C_1$ and $C_2$ of $\bar C_1
:= \{ x = 0 \}$ and $\bar C_2 := \{ y = 0 \}$ as curvettes of the
exceptional curves $E_1$ and $E_2$, respectively. Moreover, the
strict transform of $f^{-1} \{ 0 \}$ can be considered as a
curvette $C_m$ for $E_m$.

 It is well known that the multiplicities of $E_1,E_2$
and $E_m$ in $\div(\pi^\ast f)$ are $p,q$ and $pq$, respectively,
and that the multiplicity of $E_m$ in $\div(\pi^\ast dx \wedge
dy)$ is $p+q-1$. Since the matrix $(a_{ij})$ is symmetric (where
we use the notation of (2.3)) we have $a_{1m} = a_{m1} = p$ and
$a_{2m} = a_{m2} = q$. Consequently the multiplicity of $E_m$ in
$\div(\pi^\ast \omega_{ij})$ is $p + q -1 + (i-1) p + (j-1)q = ip
+ jq -1$.

The only monodromy eigenvalue on $H^0(F_0,\Bbb C)$ is $1$, and
then by A'Campo's formula the eigenvalues on $H^1(F_0, \Bbb C)$
are $$\exp \{ - 2 \pi \sqrt{-1} (\frac i q + \frac j p) \}$$ for
$1 \leq i \leq q-1$ and $1 \leq j \leq p-1.$

On the other hand, for example by [Ve2], we have already that $-1$
and $- \frac{ip + jq}{pq} = -(\frac i q + \frac jp)$ are the only
candidate poles of $Z_{top,0}(f, \omega_{ij};s)$. In fact they
really are poles which immediately implies (1) and (2).

An elegant way to check this is the formula in [Ve3, Theorem 3.3]
which yields the following compact expression (this formula
remains valid in the context of arbitrary differential forms
$\omega$):

$$Z_{top,0} (f,\omega_{ij};s) = \frac{1}{ip + jq + spq} (-1 +
\frac{1}{1+s} + \frac qi + \frac pj)$$
$$= \frac{jq + ip + (jq + ip - ij)s}{(ip + jq + spq)(1+s)}. \quad
\square$$
\enddemo

\bigskip
\noindent {\bf 2.8.} {\sl Example.} Let $f = (y^2 - x^3)^2 - x^6y$
on $(\Bbb C^2,0)$. This is one of the simplest irreducible
singularities with two Puiseux pairs. The minimal embedded
resolution of $f^{-1} \{ 0 \}$ is described in Figure 3. The
numbers $(\nu_i,N_i)$ denote as usual $1+$ (the multiplicity of
$E_i$ in $\div(\pi^\ast dx \wedge dy)$) and the multiplicity of
$E_i$ in  $\div(\pi^\ast f)$, respectively.

By A'Campo's formula the monodromy eigenvalues of $f$ at $0$ are
$1$ and all primitive roots of unity of order 6, 10, 12 and 30.
Take $\omega_{ij} = x^{i-1} y^{j-1} dx \wedge dy$ for $i,j \geq
1$. We checked that the statement in (2.5) is valid for example
for the sets of differential forms $\{ \omega_{ij} | 1 \leq i \leq
5$ and $1 \leq j \leq 3 \} \cup \{ \omega_{34} \}$ and $\{
\omega_{ij} | 1 \leq i \leq 3$ and $1 \leq j \leq 5 \} \setminus
\{ \omega_{24} \}$.

\vskip 1truecm
 \centerline{
\beginpicture
\setcoordinatesystem units <.5truecm,.5truecm>
 \putrule from -4 0 to 4 0
 \putrule from 3 -3 to 3 3
 \putrule from 2 -1.8 to 6 -1.8
 \putrule from 2 1.8  to 6 1.8
 \putrule from -3 -3 to -3 3
 \putrule from -2 -1.8 to -6 -1.8
 \putrule from -2 1.8  to -6 1.8
 \put {$E_1(2,4)$} at 7 -1.3
  \put {$E_2(3,6)$} at 7 2.3
  \put {$E_3(5,12)$} at 3 -3.6
  \put {$E_4(6,14)$} at 0 .5
 \put {$E_5(7,15)$} at -7 -1.3
  \put {$E_0(1,1)$} at -7 2.3
  \put {$E_6(13,30)$} at -3 -3.6
\endpicture}
 \vskip 1truecm
  \centerline{\smc Figure 3}
 \vskip 1truecm

\bigskip
\noindent {\bf 3. Arbitrary dimension}

\bigskip
\noindent {\bf 3.1.} First we construct a higher dimensional
generalization of the notion of curvette such that an analogue of
Proposition 2.3 is still valid.

Let $X_0$ be a smooth quasi-projective (complex) variety of
dimension $n$ and let
$$X_0 \overset{\pi_1} \to \longleftarrow X_1 \overset{\pi_2} \to \longleftarrow X_2
\overset{\pi_3} \to \longleftarrow \cdots \overset{\pi_i} \to
\longleftarrow X_i \overset{\pi_{i+1}} \to \longleftarrow \cdots
\overset{\pi_m} \to \longleftarrow X_m$$ be a composition $\pi$ of
$m$ blowing-ups $\pi_i$ with smooth irreducible centre
$Z_{i-1}(\subset X_{i-1})$ having normal crossings with the
exceptional locus of $\pi_1 \circ \cdots \circ \pi_{i-1}$. Denote
the exceptional locus of $\pi_i$, as well as its consecutive
strict transforms, by $E_i$.

Recall that, when created, $E_i$ has the structure of a $\Bbb
P^k$-bundle $E_i \overset{p_i} \to \longrightarrow Z_{i-1}$, where
$k = n-1- \dim Z_{i-1}$. We have $\Pic E_i \cong \Bbb Z L_i \oplus
p^\ast_i \Pic Z_{i-1}$, where $L_i$ is the divisor class
corresponding to the canonical sheaf ${\Cal O}_{E_i}(1)$ on $E_i$.
The self-intersection $E^2_i$ of $E_i$ on $X_i$, considered in
$\Pic E_i$, is equal to $-L_i$ [Ha, Theorem II 8.24]. (When
$Z_{i-1}$ is a point, $L_i$ is just the hyperplane class on $E_i
\cong \Bbb P^{n-1}$.)

\bigskip
\proclaim {\bf 3.2. Proposition} One can construct consecutively
for $j = 1,\cdots,m$ a smooth hypersurface $C_j$ on $X_j$ such
that

(1) $C_j$ has normal crossings with $E_1 \cup E_2 \cup \cdots \cup
E_j$, with (the strict transform of) previously created
$C_1,\cdots,C_{j-1}$, and with the next centre of blowing-up $Z_j$
(and such that $Z_j \not\subset C_j$);

(2) in $\Pic E_j$ we have $C_j \cap E_j = L_j + p_j^\ast B_j$ for
some $B_j \in \Pic Z_{j-1}$;

(3) denoting $\tilde C_j := \pi_j(C_j) \subset X_{j-1}$, we have
$\pi^\ast_j \tilde C_j = E_j + C_j$ in $\Pic X_j$. So the
multiplicity of $\tilde C_j$ along $Z_{j-1}$ is $1$.

Note that by (1) the strict transforms in $X_m$ of all $E_j$ and
$C_j$ form a normal crossings divisor.

(4) Given another hypersurface $H$ on $X_m$ having normal
crossings with $E_1 \cup \cdots \cup E_m$, we can choose
$C_1,\cdots,C_m$ such that furthermore $H$ and all $E_j$ and $C_j$
form a normal crossings divisor on $X_m$.
\endproclaim

\noindent \demo{Proof} Fix a $j \in \{1,\cdots,m \}$. Consider the
sheaf ${\Cal O}_{X_j}(1)$ on $X_j$, associated to the blowing-up
map $\pi_j : X_j \rightarrow X_{j-1}$. We choose an ample
invertible sheaf $\Cal L$ on $X_{j-1}$. By [Ha, II Proposition
7.10] we have for some $k > 0$ that the sheaf ${\Cal O}_{X_j}(1)
\otimes \pi^\ast_j {\Cal L}^k$ on $X_j$ is very ample over
$X_{j-1}$. So its global sections generate a base point free
linear system on $X_j$; we take $C_j$ as a general element of this
linear system. By Bertini's theorem $C_j$ satisfies (1). The
inverse image on $X_m$ of this linear system is still base point
free. So a general element will also satisfy the extra condition
in (4).

We now verify that the intersection product $C_j \cdot E_j$,
considered in $\Pic E_j$, is of the form $L_j + p^\ast_j B_j$,
which yields (2). Denote $\beta : E_j \hookrightarrow X_j$ and
$\alpha : Z_{j-1} \hookrightarrow X_{j-1}$. Since the divisor
class corresponding to ${\Cal O}_{X_j}(1)$ on $X_j$ is $-E_j$ we
have
$$\split C_j \cdot E_j & = (- E_j + \pi^\ast_j (\cdots)) \cdot E_j  = -
E^2_j + \beta^\ast \pi^\ast_j (\cdots) \\
& = L_j + p^\ast_j \alpha^\ast (\cdots). \endsplit$$

Finally we verify (3). Certainly $\pi^\ast_j \tilde C_j = \mu E_j
+ C_j$ where $\mu$ is the multiplicity of $\tilde C_j$ along
$Z_{j-1}$. Intersecting with $E_j$ yields $\beta^\ast \pi^\ast_j
\tilde C_j = \mu (-L_j) + C_j \cdot E_j$, and hence by the
previous calculation $p^\ast_j (\cdots) = - \mu L_j + L_j +
p^\ast_j (\cdots)$. So indeed $\mu = 1$. \quad $\square$
\enddemo

$$\CD
E_j @>\beta>> X_j \\
@Vp_jVV       @VV\pi_jV \\
Z_{j-1} @>\alpha>> X_{j-1}
\endCD
$$
\bigskip
\noindent {\bf 3.3.} The $C_j$ constructed above satisfy an
anologous statement as Proposition 2.3 for curves. For the proof
however we need another approach.

\bigskip
\proclaim{Proposition} We use the notation of (3.1) and (3.2).
Denote also $\bar C_i := \pi (C_i) \subset X_0$ and $\pi^\ast \bar
C_i = \sum^m_{j=1} a_{ij} E_j + C_i$ for $i = 1,\cdots,m$. Then
the determinant of the $(m \times m)$-matrix $(a_{ij})$ is equal
to 1. In particular $\gcd_{1 \leq i \leq m} \{ a_{ij} \} = 1$ for
all $j$. \endproclaim

\noindent \demo{Proof} We proceed by induction on $m$. When $m=1$
we have $\pi^\ast \bar C_1 = E_1 + C_1$ by Proposition 3.2(3), and
so indeed $a_{11} = 1$. Take now $m > 1$. By the same proposition
we have $\pi^\ast_m \tilde C_m = 1 E_m + C_m$, where $\tilde C_m =
\pi_m (C_m)$. Say $Z_{m-1}$ is contained in precisely $E_j, j \in
J( \subset \{ 1,\cdots,m-1 \})$. Then
$$\pi^\ast \bar C_m = \pi^\ast_m (\sum^{m-1}_{j=1} a_{mj} E_j +
\tilde C_m) = (\sum_{j \in J} a_{mj} + 1) E_m + (\cdots),$$ saying
that $a_{mm} = \sum_{j \in J} a_{mj} + 1$. On the other hand,
since $Z_{m-1}$ is not contained in (the strict transform of) any
$C_1,C_2,\cdots,C_{m-1}$, we have $a_{im} = \sum_{j \in J} a_{ij}$
for $i = 1,\cdots,m-1$. Hence

$$\det(a_{ij}) \Sb 1 \leq i \leq m \\ 1 \leq j \leq m \endSb = \vmatrix
a_{11} & \cdots & a_{1,m-1} & 0 \\
\vdots & \ddots & \vdots & \vdots \\
a_{m-1,1} & \cdots & a_{m-1,m-1} & 0 \\
a_{m1} & \cdots & a_{m,m-1} & 1
\endvmatrix = \det (a_{ij}) \Sb 1 \leq i \leq m-1  \\ 1 \leq j \leq m-1 \endSb = 1,$$
where the last equality is the induction hypothesis. \quad
$\square$ \enddemo

\bigskip
\noindent {\sl Note.} In contrast to the curve case, the matrix
$(a_{ij})$ is in general not symmetric in higher dimensions, even
when all $\pi_i$ are point blowing-ups. Take for example $\dim X_0
= 3$, $Z_1$  a point on $E_1$, and $Z_2$ a point on $E_1 \cap
E_2$. Then

$$(a_{ij}) = \pmatrix 1 & 1 & 2 \\ 1 & 2 & 3 \\ 1 & 2 & 4
\endpmatrix .$$

\bigskip
\noindent {\bf 3.4.} We now present higher dimensional versions of
Theorem 2.4. We first look at zeroes or poles of monodromy zeta
functions. According to Lemma 1.3, this way we treat in fact all
monodromy eigenvalues. For isolated singularities we present a
finer result.

\bigskip \proclaim{3.5. Theorem} Let $f : (\Bbb C^n,0) \rightarrow
(\Bbb C,0)$ be a non-zero polynomial function (germ). Let
$\lambda$ be a zero or a pole of the monodromy zeta function of
$f$ at $0$. Then there exists a differential $n$-form $\omega$ on
$(\Bbb C^n,0)$ such that $Z_{top,0}(f,\omega;s)$ has a pole $s_0$
satisfying $\exp(2 \pi \sqrt{-1} s_0) = \lambda$.
\endproclaim

\noindent \demo{Proof} Let $f : X_0(\subset \Bbb C^n) \rightarrow
\Bbb C$ be a relevant representative of $f$ in the sense that some
embedded resolution of $f^{-1} \{ 0 \} \subset X_0$ only has
exceptional components that intersect the inverse image of $0$.
Take such an embedded resolution $\pi : X_m \rightarrow X_0$,
which is a composition of $m$ blowing-ups as in (3.1). Slightly
abusing notation, $E_i$ can now denote an exceptional component of
$\pi$ or an irreducible component of the strict transform. As
usual $N_i$ and $\nu_i - 1$ are the multiplicities of $E_i$ in the
divisor of $\pi^\ast f$ and $\pi^\ast(dx_1 \wedge \cdots \wedge
dx_n)$, respectively.

Say $\lambda$ is a $d$-th root of unity. By A'Campo's formula
there exists an exceptional component $E_{j_0}$ with $d|N_{j_0}$
and $\chi(E^\circ_{j_0} \cap \pi^{-1} \{ 0 \}) \ne 0$.

We take for $i = 1,\cdots,m$ smooth hypersurfaces $C_i$ as in
Proposition 3.2 (considered in $X_m$); as extra hypersurface in
3.2(4) we take the strict transform of $f^{-1} \{ 0 \}$. Say the
images $\bar C_i$ in $X_0$ of the $C_i$ have (reduced) equation
$g_i = 0$. As before we denote $\pi^\ast \bar C_i = \sum^m_{j=1}
a_{ij} E_j + C_i$; the $(a_{ij})$ satisfy Proposition 3.3. We take
for the moment a differential form $\omega$ of the form
$(\prod_\ell g^{m_\ell}_\ell) dx_1 \wedge \cdots \wedge dx_n$. The
multiplicity of $E_i$ in the divisor of $\pi^\ast \omega$ is
$\nu_i - 1 + \sum_\ell a_{\ell i} m_\ell$.

The candidate pole for $Z_{top,0} (f,\omega;s)$ associated to
$E_{j_0}$ is $s_0 = - \frac{\nu_{j_0} + \sum_\ell a_{\ell j_0}
m_\ell}{N_{j_0}}$. As in the proof of Theorem 2.4 we want to find
suitabe $m_\ell$. Completely analogously as in that proof, this
time using Proposition 3.3,

(1) we find a lattice of nonnegative $m_\ell$ `mod $N_{j_0}$' such
that $s_0$ satisfies $\exp(2 \pi \sqrt{-1} s_0) = \lambda$, and

(2) for `most' such $m_\ell$ the candidate poles associated to
other $E_i$ are different from $s_0$.

\noindent The argument showing that $s_0$ is really a  pole for
suitable such $m_\ell$ is more subtle now. We introduce some
notation to describe the residue of $s_0$.

\bigskip
Let $C_\ell, \ell \in J_0$, be the hypersurfaces $C_i$ that
intersect $E_{j_0} \cap \pi^{-1} \{ 0 \}$ (in $X_m$). Denote
$C^0_J := ((\cap_{j \in J} C_j) \setminus (\cup_{i \in J_0
\setminus J} C_i)) \cap (E^\circ_{j_0} \cap \pi^{-1} \{ 0 \})$ for
$J \subset J_0$; in particular $C^\circ_\emptyset = (E^\circ_{j_0}
\cap \pi^{-1} \{ 0 \}) \setminus \cup_{i \in J_0} C_i$. These
$C^\circ_J$ form a locally closed stratification of $E^\circ_{j_0}
\cap \pi^{-1} \{ 0 \}$. The residue of $s_0$ is of the form
$$\frac{1}{N_{j_0}} \biggl( \chi (C^\circ_\emptyset) + \sum_{\emptyset \ne J \subset
J_0} \chi (C^\circ_J) \prod_{j \in J} \frac{1}{1 + m_j} + \text
{contribution  of } (E_{j_0} \setminus E^\circ_{j_0}) \cap
\pi^{-1} \{ 0 \} \biggr)\, ,$$
 the last contribution also being a
rational function in the $m_\ell$ of negative degree. As in the
proof of Theorem 2.4, if $\chi (C^\circ_\emptyset) \ne 0$, this
expression is never identically zero as function in the $m_\ell$,
and so non-zero for `most' choices of $(m_\ell)_\ell$.

\bigskip
We don't see how to exclude the theoretical possibility that this
expression {\sl is} identically zero (with then necessarily
$\chi(C^\circ_\emptyset) = 0$). In this case we will adapt our
choice of $\omega$ to be sure to have the desired pole.

For each $\ell$ in $J_0$ we construct, as in Proposition 3.2, not
just one hypersurface $C_\ell$, but several ones $C_{\ell 1},
C_{\ell 2}, \cdots , C_{\ell t}$, all general enough elements in
the linear system that was considered there. Then still all $E_i,
C_{\ell j}$, and other $C_\ell$ will form a normal crossings
divisor on $X_m$. Say $g_{\ell k}=0$ is the equation of the image
of $C_{\ell k}$ in $X_0$. Now we take $\omega$ of the form
$\prod_{\ell \in J_0} (\prod^t_{k=1} g^{m_{\ell k}}_{\ell k})
\prod_{\ell \notin J_0} g^{m_\ell}_\ell dx_1 \wedge \cdots \wedge
dx_n$ such that for $\ell \in J_0$ the sum $\sum^t_{k=1} m_{\ell
k}$ is an allowed $m_\ell$ `mod $N_{j_0}$' as before. The
candidate pole $s_0$ for $Z_{top,0}(f,\omega;s)$ associated to
$E_{j_0}$ is as above, it still satisfies $\exp(2 \pi \sqrt{-1}
s_0) = \lambda$, and for `most' such $m_{\ell k}$ and $m_\ell$ the
candidate poles associated to other $E_i$ are different from
$s_0$. We now verify that for some $t$ the expression for the
residue of $s_0$ is not identically zero as function in the
$m_{\ell k}$ and $m_\ell$.

Denote $L_k := (\cup_{\ell \in J_0} C_{\ell k}) \cap
(E^\circ_{j_0} \cap \pi^{-1} \{ 0 \})$ for $k = 1,\cdots,t$. Since
all $C_{\ell k}$ all {\sl general} elements we have that all
$\chi(L_k)$ are equal, that also the $\chi(L_k \cap L_{k^\prime})$
are equal for all $k < k^\prime$, and more generally that the
$\chi(L_{k_1} \cap L_{k_2} \cap \cdots \cap L_{k_s})$ are equal
for all $1 \leq k_1 < k_2 < \cdots < k_s \leq t$. The residue of
$s_0$ is of the form
$$\frac{1}{N_{j_0}} \bigl( \chi((E^\circ_{j_0} \cap \pi^{-1} \{ 0 \}) \setminus
\cup^t_{k=1} L_k) + \cdots \bigr) \, ,$$ where the other terms
form a rational function of negative degree in the $m_{\ell k}$
and the $m_\ell$. If $\chi((E^\circ_{j_0} \cap \pi^{-1} \{ 0 \})
\setminus \cup^t_{k=1} L_k) \ne 0$, then this residue is not
identically zero and we are done. Finally we show that this must
be the case for some $t$.

Because of the normal crossings property we have $\cap^T_{k=1} L_k
= \emptyset$ for some $T(\leq n)$. Suppose that
$\chi((E^\circ_{j_0} \cap \pi^{-1} \{ 0 \}) \setminus \cup^t_{k=1}
L_k) = 0$ for {\sl all} $t = 1,\cdots,T$. These $T$ conditions can
be rewritten as
$$\cases \chi(E^\circ_{j_0} \cap \pi^{-1} \{ 0 \}) - \chi(L_1) = 0
\\ \chi(E^\circ_{j_0} \cap \pi^{-1} \{ 0 \}) - 2 \chi(L_1) + \chi (L_1 \cap L_2) =
0 \\
\chi(E^\circ_{j_0} \cap \pi^{-1} \{ 0 \}) - 3 \chi(L_1) + 3 \chi
(L_1 \cap L_2) - \chi(L_1 \cap L_2 \cap L_3) = 0 \\
\cdots \\
\chi(E^\circ_{j_0} \cap \pi^{-1} \{ 0 \}) - (T-1) \chi (L_1) +
\cdots + (-1)^{T-1} \chi (\cap^{T-1}_{k=1} L_k) =
0 \\
\chi (E^\circ_{j_0} \cap \pi^{-1} \{ 0 \}) - T \chi (L_1) + \cdots
+ (-1)^{T-1} T \chi (\cap^{T-1}_{k=1} L_k) + (-1)^T \cdot 0 = 0 \,
.
\endcases
$$
One easily verifies that the $(T \times T)$-determinant of
coefficients of this homogeneous linear system of equations in the
$\chi(\cdots)$ is nonzero. Hence in particular we should have
$\chi(E^\circ_{j_0} \cap \pi^{-1} \{ 0 \}) = 0$, contradicting our
choice of $E_{j_0}$. \quad $\square$
\enddemo

\bigskip
\noindent \proclaim{3.6. Theorem} Let $f : (\Bbb C^n,0)
\rightarrow (\Bbb C,0)$ be a nonzero polynomial function (germ).

(a) Let $\lambda$ be a monodromy eigenvalue of $f$ at $0$. Then
there exists a differential $n$-form $\omega$ and a point $P$ in a
neighbourhood of $0$ such that $Z_{top,P}(f,\omega;s)$ has a pole
$s_0$ satisfying $exp(2 \pi \sqrt{-1} s_0) = \lambda$. Moreover,
$P$ can be chosen as a generic point in the set $\Sigma$ that was
introduced after Lemma 1.3.

(a') If the eigenvalue $\lambda$ appears {\sl only} at $0$, then
there exists a differential $n$-form $\omega$ such that
$Z_{top,0}(f,\omega;s)$ has a pole $s_0$ satisfying $exp(2 \pi
\sqrt{-1} s_0) = \lambda$.

(b) Suppose that $f^{-1} \{ 0 \}$ has an isolated singularity at
$0$, and let $\lambda$ be a monodromy eigenvalue of $f$ at $0$.
Then there exists a differential $n$-form $\omega$ such that
$Z_{top,0}(f,\omega;s)$ has a pole $s_0$ satisfying $exp(2 \pi
\sqrt{-1} s_0) = \lambda$.
\endproclaim

\noindent \demo{Proof} Parts (a) and (a') follow immediately from
Theorem 3.5 and Lemma 1.3.

Part (b) is a special case of (a') for $\lambda \neq 1$. It could
however happen in (b) that $\lambda=1$ is not a zero or a pole of
the monodromy zeta function of $f$ at $0$ (when $n$ is even). In
that case we pick any exceptional $E_{j_0}$ with
$\chi(E^\circ_{j_0}) \neq 0$, and we proceed as in the proof of
Theorem 3.5 to construct a suitable $\omega$ and a pole $s_0$ of
$Z_{top,0}(f,\omega;s)$ with $exp(2 \pi \sqrt{-1} s_0) = 1$. (Note
that the constructed $s_0$ is in this case indeed an integer.) By
A'Campo's formula there is always such an exceptional $E_{j_0}$,
except when $n$ is even and the characteristic polynomial
$P_{n-1}(t) = t-1$. E.g. by [AGLV, page 70] this implies that $f$
has a so-called non-degenerate or Morse singularity at $0$ (i.e.,
$f$ is in local analytic coordinates of the form
$y_1^2+y_2^2+\dots+y_n^2$). In this easy special case one has
$$Z_{top,0}(f,dx_1\wedge \cdots \wedge dx_n;s) = \frac{n}{(1+s)(n+2s)}.\quad \square$$
\enddemo

\noindent {\bf 3.7.} As in the curve case the zeta functions
constructed in the proof of Theorem 3.5 can in general have other
poles that don't induce monodromy eigenvalues of $f$, and it would
be interesting to study in arbitrary dimension the validity of the
`principle' in Note 2.5. We present an example below.

\bigskip
\noindent {\bf 3.8.} {\sl Example.} Let $f = x^d + y^d + z^d$ on
$(\Bbb C^3,0)$ with $d \geq 3$. Blowing up the origin yields an
embedded resolution $\pi$ of $f^{-1} \{ 0 \}$. The exceptional
surface $E \cong \Bbb P^2$ has multiplicities 2 and $d$ in
$\div(\pi^\ast dx \wedge dy \wedge dz)$ and $\pi^\ast f$,
respectively. It intersects the strict transform in a smooth curve
$D$ of degree $d$ and hence with Euler characteristic $3d-d^2$. By
A'Campo's formula the monodromy zeta function of $f$ at $0$ is
$$\zeta_{f,0}(t) = (t^d-1)^{d^2-3d+3};$$
and the monodromy eigenvalues of $f$ at $0$ are precisely all
$d$-th roots of unity. Take $\omega_i := x^{i-1} dx \wedge dy
\wedge dz$ for $1 \leq i \leq d$. The strict transform of $\{ x =
0 \}$ intersects $E$ in a line; this line intersects $D$
transversely in $d$ points. Hence
$$\split
Z_{top,0}(f,\omega_i;s) & = \frac{1}{(2+i)+sd} \biggl( (d-1)^2 +
\frac{2-d}{i} + \frac{2d-d^2}{1+s} + \frac{d}{i(1+s)} \biggr) \\
& = \frac{(i(d-1)^2 + 2-d)s + (2+i)}{i((2+i)+sd)(1+s)}.
\endsplit
$$
When $i \ne d-2$, one easily verifies that the two candidate poles
$-1$ and $- \frac{2+i}{d}$ are really poles. When $i = d-2$ we
have
$$Z_{top,0} (f,\omega_{d-2};s) = \frac{1 +
(d-2)^2s}{(d-2)(1+s)^2}$$ and $-1$ is a pole (of order 2 if $d >
3$ and of order 1 if $d=3$).

So the set of differential forms $\{ \omega_i|1 \leq i \leq d \}$
satisfies the analogous principle as in (2.5). We can even delete
$\omega_{d-2}$ from this set.

\enddocument